\documentclass[12pt]{article}
\usepackage[tbtags]{amsmath}
\usepackage[]{amssymb}
\usepackage[]{amsbsy}
\usepackage[]{amsthm}
\usepackage[dvips]{graphicx}

\newtheorem{theorem}{Theorem}[section]

\newtheorem{lemma} {Lemma} [section]
\newtheorem{defi}{Definition} [section]
\newtheorem{example}{Example} [section]
\newtheorem{remark}{Remark} [section]

\def\<{\langle}
\def\>{\rangle}

\def\be{\begin{equation}}
\def\bea{\begin{eqnarray}}

\def\ee{\end{equation}}
\def\eea{\end{eqnarray}}

\def\l{\lambda}

\def\rem#1{\begin{remark}{\rm #1}\end{remark}}

\def\splt#1{\begin{equation} \begin{split} #1 \end{split}\end{equation}}

\renewcommand{\(}{\left(}
\renewcommand{\)}{\right)}

\begin{document}
\date{}
\baselineskip=1.2\normalbaselineskip
\thispagestyle{empty}

\title{{\noindent \bf  Boundary Feedback Control of Unstable Heat Equation with
Space and Time Dependent Coefficient }}

\author{
{\sc   Junji Jia  } \\
 \normalsize{ Department of Applied Mathematics}\\
 \normalsize{ Room 210, Middlesex College, The University of Western Ontario }\\
  \normalsize{London, Ontario, N6A 5B7, Canada }\\
 \normalsize{\tt Email: jjia5@uwo.ca} \vspace*{-2ex}}
\maketitle

\vspace*{2ex}
\begin{abstract}
In this paper, we study the control of the linear heat equation
with a space and time dependent coefficient function
by the Dirichlet and Neumann boundary control laws. This
equation models the heat diffusion and space, time dependent heat generation in
a one dimensional rod. Without control, the system is unstable if the
coefficient function is positive and large. With boundary control based on the
state feedback, we show that for the time analytic coefficient
$a(x,t)$, the exponential stability of the system at any rate
can be achieved. It is further shown that both the control of
the Dirichlet and Neumann boundary value systems can be stabilized using this
method. In doing this, the control kernels are {\it
explicitly} calculated as series of approximation and they are used in the
simulations. The numerical simulation confirmed the theoretical arguments and
the controllability of the system. The possible future works are discussed at
the end. This is the continuation
of the recent work of Liu [{\it SIAM J. Cont. Optim.} vol. 42, pp. 1033-1043] and Smyshlyaev and Krstic [accepted by
{\it Automatica}].
\end{abstract}



\thispagestyle{empty}
\newpage
\section{Introduction}
\hspace{15pt}
In this paper, we continue the study of feedback control of an unstable heat
equation \be u_t(x,t)=u_{xx}(x,t)+au(x,t) \mbox {  in
$(0,1)$.}  \label{gene}\end{equation}
This equation arises from the heat generation and diffusion of a rod in which
the heat generation is only only sensitive to the position $x$ but also the time
$t$. A simple version of this
feedback control problem where $a=\mu$, a constant, was first addressed by
Boskovic, Krstic and Liu in \cite{bkl} and the instability when $\mu$ is large
was shown there. Later, Balogh and Krstic \cite{bk0, bk} studied the case where
$a=a(x)$, a function dependent on $x$. Since they used a backstepping
approximation method, the same problem was reconsidered and solved by Liu
\cite{liunew} using a transformation skill. More recently, the generalization
to $a=a(t)$ was explicitly considered by Smyshlyaev and Krstic \cite{am2}. In
our paper, we generalize it to the case of $a=a(x,t)$.

For the theories that use boundary feedback control to stabilize the parabolic
equations, one can refer to \cite{bkl} and \cite{bk} and references there. In
comparison of the existing literature, the novelty of this paper is the
"explicit" construction of the feedback laws of $a=a(x,t)$ case for some
$a(x,t)$ and the proof
that both the kernel for the Dirichlet and Neumann boundary control problem
exist.

The paper is organized as follows. Section \ref{dir} is devoted to
the stabilization of Dirichlet boundary value problems of $a=a(x,t)$ case,
Section \ref{neum} to the stabilization of Neumann boundary value problems. In
Section \ref{simu}, the simulation of $a(x,t)=x(bt+c)$ is studied and the
validity of the control are shown. In section 5, we give some remarks of
interesting phenomena, which may be helpful for the future works.

In what follows, we denote the region $(0,1)_x$ by $X$, $(0,1)_t$ by $T$,
$0\leq y\leq x\leq 1$ by $\Omega$ and
$X\times T$ by $XT$ and $\Omega\times T$ by $\Omega T$ and so on. We also denote by
$H^s\! :\!XT$ the usual Sobolev space (see, e.g.,  \cite{Ada, LM})
for  any $s\in {\bf R}$.  For $s\geq 0$, $H^s_0\! :\!XT$   denotes the completion of
$C_0^\infty\! :\!XT$ in $H^s\! :\!XT$, where $C_0^\infty \! :\!XT$ denotes the space of all
infinitely differentiable functions on $XT$ with compact support  in $XT$.
 We denote by $\|\cdot\|$   the norm  of
$L^2\! :\!XT$ with respect to space variable $x$. $C^n\! :\!XT$ denotes the space of all
$n$ times  continuously differentiable functions on $XT$.

\section{Dirichlet Boundary Conditions}\label{dir}
\setcounter {equation} {0}
\setcounter {theorem} {0}
\setcounter {lemma} {0}

It is well known that the Dirichlet boundary value problem
\begin{equation}\label{linear}
\left\{\begin{array} {ll}
 u_t(x,t)=  u_{xx}(x,t)+a(x,t) u(x,t)
  & \mbox { in $XT,$}\\
u(0,t)=u(1,t)=  0  & \mbox { in $T$}
\end{array}
\right.
\end{equation}
is unstable if function $a$ is positive and large \cite{liunew}.

To design a boundary feedback law to stabilize the above system for
function $a(x,t)$, we consider the problem
\be\label{liupb}
\left\{
\begin{array}{ll}
k_{xx}(x, y, t)-k_{yy}(x, y, t)-k_{t}(x, y, t)=(a(y, t)+\l)k(x, y, t),
&\mbox{in $\Omega T$}\\
k(x, 0,t)=0,& \mbox{in $XT$}\\
k_x(x, x, t)+k_y(x, x, t)+\frac{d}{d x}(k(x, x, t))=
a(x,t)+\l,& \mbox{in $XT$}
\end{array}
\right.
\end{equation}
where $\l$ is any positive constant. Why we consider this problem will becomes
clear in Lemma~\ref{lem2}. For the moment, let
us assume this problem has a unique solution $k$ for some $a$. Using the solution $k$,
we then obtain a Dirichlet boundary feedback law
\be\label{dirfd}
u(1,t)= - \int_0^1 k(1,y,t)u(y, t)dy \qquad \mbox { in $T$}
\end{equation}
and Neumann boundary feedback law
\be\label{neumfd}
u_x(1,t)= -k(1,1)u(1,t)- \int_0^1 k_x(1,y,t)u(y, t)dy
\qquad \mbox { in $T$.}
\end{equation}
With one of the boundary feedback laws, the  system
\begin{equation}
\left\{\begin{array} {ll}
 u_t(x,t)=  u_{xx}(x,t)+a(x,t) u(x,t)
  & \mbox { in $XT,$}\\
u(0,t)=0& \mbox { in $T,$}\\
  u(x, 0)=u^0(x)     & \mbox {  in $T$ }
\end{array}
\right.
\label{dche}
\end{equation}
is exponentially stable.
In this controlled system, the left hand end of a rod is insulated while the temperature or the heat flux at the other
end is adjusted according to the measurement of $k$-weighted
averaged temperature over the whole rod. Physically, if the
destabilizing heat is generating inside the rod, then we cool
the right end of the rod so that it is not overheated.
 To state this result, we also introduce the
compatible conditions for the initial data
\bea\label{compat1}
u^0(0)&=&0,\quad u^0(1,t)=- \int_0^1 k(1,y,t)u^0(y)dy\\
u^0(0)&=&0,\quad u_x^0(1,t)= -k(1,1)u^0(1,t)- \int_0^1 k_x(1,y,t)u^0(y)dy
\label{compat2}
\end{eqnarray}

\begin{theorem} \label{cstabthm2} Assume that  $\l >0$ is any positive
constant and $a(x,t)$ is continuous in
$x$ and analytic in $t$ uniformly in $x$.
For arbitrary  initial  data $u^0(x) \in H^1X $ with
compatible condition \eqref{compat1} or \eqref{compat2}, equation \eqref{dche} with either \eqref{dirfd} or \eqref{neumfd} has a unique solution that satisfies
 \be\label{cstab4}
\|u(t)\|_{H^1} \leq M \|u^0\|_{H^1}e^{- \l  t }, \quad t\in(0,t_0),
t_0<1,
\end{equation}
where $M$ is a positive constant independent of $u^0$.
\end{theorem}

The idea of proving the theorem is to carefully construct a transformation
\[
w(x,t)=u(x,t)+\int_0^x k(x, y, t)u(y,t)dy
\]
 to convert the system \eqref{dche} with either \eqref{dirfd} or \eqref{neumfd} into
the exponentially stable system
\begin{equation}
\left\{\begin{array} {ll}
 w_t= w_{xx}-\l w
  & \mbox { in $XT,$}\\
  w(0,t) = w(1,t)=0& \mbox { in $T,$}\\
  w(x, 0)=w^0(x) & \mbox {  in $X,$ }
\end{array}
\right.
\label{stab}
\end{equation}
or
\begin{equation}
\left\{\begin{array} {ll}
 w_t= w_{xx}-\l w
  & \mbox { in $XT,$}\\
  w(0,t) = w_x(1,t)=0& \mbox { in $T,$}\\
  w(x, 0)=w^0(x) & \mbox {  in $X,$ }
\end{array}
\right.
\label{neumstab}
\end{equation}
where $w^0(x)=u^0(x)+\int_0^x k(x, y, t)u^0(y)dy$. This will be done by the
following lemmas. The first one of these lemmas is due to Colton \cite{col,sd}. We
conveniently rewrite its statement and quote the proof here for reference.

\begin{lemma} \label{lem1} {\bf (Colton)} Suppose that $a(x,t)$ is continuous in
$x$ and analytic in $t$ uniformly in $x$. Then
problem \eqref{liupb} has a unique solution which is twice continuously
differentiable in $x,~y$.
\end{lemma}

\begin{proof}
Using the substitutions
\[x=\xi+\eta,~y=\xi-\eta,~a(y,t)=a(\xi-\eta,t),~G(\xi,\eta,t)=k(x,y,t),\]
the problem (\ref{liupb}) is transformed to
\begin{equation} \left\{ \begin{array}{ll}
G_{\xi\eta}(\xi,\eta,t)=\left(a\left(\xi-\eta,t\right)+\lambda+\frac{\partial}{\partial
t}\right)G(\xi,\eta,t)&0\leq \eta\leq \xi\leq 1,t\in T\\
G(\xi,\xi,t)=0&0\leq \xi\leq 1,t\in T\\
\frac{\partial}{\partial
\xi}(G(\xi,0,t)=\frac{1}{2}\left(a\left(\xi,t\right)+\lambda\right)& 0\leq
\xi\leq 1,t\in T.\end{array} \right. \nonumber\end{equation}
Integrating over $\xi$ and $\eta$, we can find
\begin{equation}
G(\xi,\eta,t)=\frac{1}{2}\int_{\eta}^{\xi}\left(a(\tau,t)+\lambda\right)d\tau+\int_{\eta}^{\xi}\int_0^{\eta}\left(
a\left(\tau-s,t\right)+\lambda+\frac{\partial}{\partial
t}\right)G(\tau,s,t)dsd\tau~.\label{eqofg}\end{equation}

Then Colton showed that the above equation has
a solution using the ``method of dominants''. This method works as follows.
If we are  given two series
\be S_1=\sum_{n=1}^{\infty}a_{1n}t^n,~S_2=\sum_{n=1}^{\infty}a_{2n}t^n,~t\in
(0,1),\ee
where $a_{2n}\geq0$, then we say $S_2$ {\it dominates} $S_1$ if $|a_{1n}|\leq
a_{2n},~n=1,2,3,\cdots$, and write $S_1\ll S_2$. It can be easily checked that
\bea
\mbox{if } S_1\ll S_2,&&\mbox{then }|S_1|\leq S_2,\label{eq:abspro}\\
&&\mbox{ and }\frac{\partial S_1}{\partial t}\ll\frac{\partial S_2}{\partial t},~S_1\ll S_2(1-t)^{-1};\\
\mbox{if } S_1\ll S_2,S_2\ll S_3,&&\mbox{then }S_1\ll S_3;\\
\mbox{if } S_1\ll S_2,S_3\ll S_4,&&\mbox{then }S_1+S_2\ll S_3+S_4.\eea
Moreover, the property of ``dominant'' can also be kept if the integrals of the
two series are not with respect to $t$ but other variables, that is
\be \mbox{if } S_1\ll S_2,\mbox{ then
}\int_a^bS_1dx\ll\int_a^bS_2dx.\label{eq:intpro}\ee
On the other hand, if a function $f$ is analytic in $t\in(0,1)$, then there
exist a positive constant $C$ such that $f\ll C(1-t)^{-1}$.

Using this method in our problem, it can be shown that equation (\ref{eqofg}) has a
unique twice continuously differentiable solution if $a$ is analytic in $t$. In
fact, since $a(x,t)$ is analytic, we can let $C$ be a positive constant such
that we have
\begin{equation}a(\xi-\eta,t)+\lambda\ll C(1-t)^{-1}~.\nonumber\end{equation} By the
properties of dominant (\ref{eq:abspro})-(\ref{eq:intpro}), it can be shown by
induction that for the following series
\begin{equation}
G(\xi,\eta,t)=\sum_{n=0}^{\infty}G_n(\xi,\eta,t)\label{series}\end{equation}
where
\begin{eqnarray*}
G_0(\xi,\eta,t)&=&\frac{1}{2}\int_{\eta}^{\xi}\left(a(\tau,t)+\lambda\right)d\tau\\
G_{n+1}(\xi,\eta,t)&=&\int_{\eta}^{\xi}\int_0^{\eta}\left(
a\left(\tau-s,t\right)+\lambda+\frac{\partial}{\partial
t}\right)G_n(\tau,s,t)dsd\tau~,\end{eqnarray*}
we have
\begin{equation}\nonumber
G_n\ll\frac{2^nC^n\xi^n\eta^n}{n!}\left(1-t\right)^{-(n+1)}\end{equation}
and hence
\begin{equation}
|G_n|\le\frac{2^nC^n\xi^n\eta^n}{n!}\left(1-t\right)^{-(n+1)}~.\label{eq:gn}\end{equation}
Thus the series (\ref{series}) converges absolutely and uniformly. Moreover,
since $a$ is $C^1$, $G$ is clearly $C^2$ in $x,~y$.

\end{proof}

\rem{The proof of Lemma~\ref{lem1} provides a numeric computation scheme of successive approximation to compute
the kernel function $k$ in our feedback laws \eqref{dirfd}
and \eqref{neumfd}. This makes the feedback laws \eqref{dirfd}
and \eqref{neumfd} very useful in real problems. }

\begin{lemma} \label{lem2} Let $k(x, y, t)$ be the solution of
problem \eqref{liupb} and define the linear bounded operator $K: H^i\! :\!XT_0\to
H^i\! :\!XT_0$ ($i=0,1,2$) where $T_0$ denotes $(0,t_0)(t_0<1)$ by
\be \label{trans}
w(x, t)=(Ku)(x,t)=u(x, t)+\int_0^x k(x, y, t)u(y, t)dy, \quad \mbox{ for
$u\in H^i\! :\!XT_0$}.
\end{equation}
Then,
\begin{enumerate}
\item $K$ has a linear bounded inverse $K^{-1}:~ H^i\! :\!XT_0\to H^i\! :\!XT_0$ ($i=0,1,2$).
\item $K$ converts the system \eqref{dirfd} and \eqref{dche} and system \eqref{neumfd} and \eqref{dche} into
\eqref{stab} and \eqref{neumstab}, respectively.
\end{enumerate}
\end{lemma}

\begin{proof}
To prove {\it 1.}: that  \eqref{trans} has a bounded inverse, we set
\[
v(x, t)=\int_0^x k(x, y, t)u(y, t)dy
\]
and then
\[
w(x, t)=u(x, t)+v(x, t).
\]
Hence we have
\splt{\label{inveq}
v(x, t)&=\int_0^x k(x, y, t)[w(y, t)-v(y, t)]dy\\
&=\int_0^x k(x, y, t)w(y, t)dy-\int_0^x k(x, y, t)v(y, t)dy.
}
To show that
this equation has a unique continuous solution, we set
\begin{eqnarray*}
v_0(x,t)&=& \int_0^x k(x, y, t)w(y, t)dy,\\
v_n(x,t)&=&-\int_0^x k(x, y, t)v_{n-1}(y, t)dy
\end{eqnarray*}
Though from (\ref{eq:gn}) we see that $k(x,y,t)$ is divergent with $t$
approaches $1^-$, it
is still clear that under the requirement $t\in (0,t_0)$, the $k(x,y,t)$ is
still bounded in $\Omega T_0$. Thus we can denote $M=\displaystyle\sup_{0\le y\le x\le 1, t\in
T_0} |k(x, y, t)|$.
Then one can see
\begin{eqnarray*}
|v_0(x,t)|&\le& M\|w\|,\\
|v_1(x,t)|&\le& M^2\|w\| x,\\
|v_2(x,t)|&\le& \frac{M^3\|w\|}{2!}x^2,
\end{eqnarray*}
and by induction,
\[
|v_n(x,t)|\le \frac{M^{n+1}\|w\|}{n!}x^{n}.
\]
These estimates show that the series
\[
v(x, t)=\sum_{n=0}^\infty v_n(x,t)
\]
converges absolutely and uniformly in $XT_0$
and that its sum is a continuous solution of equation \eqref{inveq}. Moreover,
there exists a constant $C>0$ such that
\be\label{l2bd}
\|v\|\le C \|w\|.
\end{equation}
This implies that there exists a bounded linear operator
 $\Phi: L^2\! :\!XT_0\to L^2\! :\!XT_0$ such that
\[
v(x, t) = (\Phi w)(x, t)
\]
and then
\be\label{inv}
u(x, t) = w(x, t)-v(x, t)=((I-\Phi) w)(x, t)=(K^{-1}w)(x, t).
\end{equation}
It is clear that $K^{-1}: L^2\! :\!XT_0\to L^2\! :\!XT_0$ is bounded.
To show that  $K^{-1}: H^1\! :\!XT_0\to H^1\! :\!XT_0$ is bounded, we
take derivative in \eqref{inveq} and obtain
\splt{\nonumber
&\quad v_x(x, t)\\
&=k(x, x, t)w(x, t)+\int_0^x k_x(x, y, t)w(y, t)dy-k(x, x, t)v(x, t)-\int_0^x k_x(x, y, t)v(y, t)dy,
}
which, combining \eqref{l2bd}, implies that there exists  constant $C>0$ such that
\[
\|v_x\|\le C\|w\|
\]
and then by \eqref{inv}
\[
\|u\|_{H^1}\le \|w\|_{H^1}+\|v\|_{H^1}\le C\|w\|_{H^1}.
\]
By analogy, we can show that  $K^{-1}: H^2\! :\!XT_0\to H^2\! :\!XT_0$ is bounded.

To prove {\it 2.}: that the transformation \eqref{trans} converts the system \eqref{dirfd} and \eqref{dche} and system \eqref{neumfd} and \eqref{dche} into
\eqref{stab} and \eqref{neumstab}, respectively, we compute as follows.
\splt{\label{comp1}
w_t(x,t)&=u_t(x,t)+\int_0^x k(x, y, t)u_t(y,t)dy+\int_0^x k_t(x, y, t)u(y,t)dy\\
&=u_t(x,t)+\int_0^x k(x, y, t)[u_{yy}(y,t)+a(y,t) u(y,t)]dy+\int_0^x k_t(x, y, t)u(y,t)dy\\
&=u_t(x,t)+k(x, x, t)u_{x}(x,t)-k(x, 0,t)u_{x}(0,t)\\
&\quad-
k_y(x, x, t)u(x,t)+k_y(x, 0,t)u(0,t)\\
&\quad+\int_0^x [k_{yy}(x, y, t)+k_t(x,y,t)+k(x, y, t)a(y,t)] u(y,t)dy,
}
\splt{\label{comp2}
w_x(x,t)&=u_x(x,t)+k(x, x, t)u(x,t)+\int_0^x k_x(x, y, t)u(y,t)dy,
}

\splt{\label{comp3}
w_{xx}(x,t)&=u_{xx}(x,t)+\frac{d}{d x}(k(x, x, t))u(x,t)+k(x, x, t))u_x(x,t)\\
&\quad+k_x(x, x, t)u(x,t)+\int_0^x k_{xx}(x, y, t)u(y,t)dy.
}
It then follows from  \eqref{liupb} and \eqref{dche} that
\splt{\label{comp4}
&\quad w_t-w_{xx}+\l w\\
&=u_t(x,t)+k(x, x, t)u_{x}(x,t)-k(x, 0,t )u_{x}(0,t)\\
&\quad-k_y(x, x, t)u(x,t)+k_y(x, 0,t)u(0,t)\\
&\quad+\int_0^x [k_{yy}(x, y, t)+k_t(x,y,t)+k(x, y, t)a(y,t)] u(y,t)dy\\
&\quad -u_{xx}(x,t)-\frac{d}{d x}(k(x, x, t))u(x,t)-k(x, x, t))u_x(x,t)\\
&\quad-k_x(x, x, t)u(x,t)-\int_0^x k_{xx}(x, y, t)u(y,t)dy\\
&\quad +\l u(x,t)+\l \int_0^x k(x, y, t)u(y,t)dy\\
&=\left( a(x,t)-k_x(x, x, t)-k_y(x, x, t)-\frac{d}{d x}(k(x, x, t))+\l\right)u(x,t)\\
&\quad+  k_y(x, 0,t)u(0,t)-k(x, 0,t)u_{x}(0,t)\\
&\quad+ \int_0^x \left[ k_{yy}(x, y, t)+k_t(x,y,t)-k_{xx}(x, y, t)+(a(y,t)+\l)k(x, y, t)\right] u(y,t)dy\\
&=0.
}
By the boundary condition of \eqref{dche}, we deduce that
$w(0,t)=0$.
Using feedback law \eqref{dirfd} or \eqref{neumfd}, we obtain
\[
w(1,t)=u(1,t)+\int_0^1 k(1,y,t)u(y,t)dy=0,
\]
or
\[
w_x(1,t)=u_x(1,t)+k(1,1)u(1,t)+\int_0^1 k_x(1,y,t)u(y, t)dy
=0.
\]
\end{proof}
\begin{remark}
We require $t\in (0,t_0)(t_0<1)$ in this lemma since $k$ is not bounded when
$t$ approaches $1^-$. However, this requirement is tolerable since
in practice, we only require the system to be stable in a time interval
$t\in(0,T_0)$ and this can be rescaled to $(0,t_0)$ by $t\to t/T(T>T_0)$ in the heat
equation. Moreover, $T_0$ can be choose as big as we want.
\end{remark}
Now we are ready to prove Theorem~\ref{cstabthm2}.

\begin{proof}[Proof of Theorem~\ref{cstabthm2}]
We first note that problem \eqref{dche} with either
\eqref{dirfd} or \eqref{neumfd} is well posed since,
by Lemma \ref{lem2}, they can be transformed to the problem \eqref{stab} or
\eqref{neumstab} via the isomorphism defined by \eqref{trans} and the problem \eqref{stab} or \eqref{neumstab} is
well posed (see, e.g., \cite[Chap. IV]{la}). Moreover, there exists a positive constant $C >0$ such that
\begin{eqnarray*}
 \|u(t)\|_{H^1}&\leq& C\|w(t)\|_{H^1},\\
 \|w^0\|_{H^1}&\leq& C\|u^0\|_{H^1}.
\end{eqnarray*}
Therefore, it is sufficient to prove \eqref{cstab4} for the solution $w$ of \eqref{stab} or
\eqref{neumstab}. We do so only for problem \eqref{stab}
since the situation for problem \eqref{neumstab} is similar.

We define the energy
  \[
 E(t)= \frac{1}{2}\int_0^1 w(x,t)^2 dx.
 \]
Multiplying the first equation of \eqref{stab} by $w$ and integrating from $0$ to $1$ by parts we get
\bea
 \dot  E(t) &=& w_x w\Big|_0^1- \int_0^1 w_x(x,t)^2 dx
- \l\int_0^1   w(x,t) ^2 dx\nonumber\\
&=& -\int_0^1 w_x(x,t)^2 dx
- \l \int_0^1   w(x,t) ^2 dx\nonumber\\
&\le& -2\l E(t),\nonumber
\end{eqnarray}
which implies
\[
 E(t)\leq E(0)e^{-2\l t} ,\;\;\;  \mbox {  for $  t\geq 0.$ }
\]
Set
\[
V(t)=\int_0^1 w_x(x,t)^2 dx.
\]
Multiplying the first equation of (\ref{stab}) by $w_{xx}$ and integrating from 0 to 1 by parts we obtain
\bea
\dot V (t) &=&  -2 \int_0^1  w_{xx} ^2 dx +2\l \int_0^1  w w_{xx }  dx\nonumber\\
 &= & -2 \int_0^1  w_{xx} ^2 dx - 2\l \int_0^1  w_{x} ^2 dx \nonumber\\
& \leq&  -2 \l V(t) ,\nonumber
\end{eqnarray}
which implies that
\[
V (t)\leq V(0)e^{ -2\l t}.
\]
This shows that \eqref{cstab4} holds.
\end{proof}

\section{Neumann Boundary Conditions} \label{neum}
\setcounter {equation} {0}
\setcounter {theorem} {0}
\setcounter {lemma} {0}

To stabilize the Neumann boundary value problem
\begin{equation}\label{neumanpro}
\left\{\begin{array} {ll}
 u_t(x,t)=  u_{xx}(x,t)+a(x,t) u(x,t)
  & \mbox { in $XT,$}\\
u_x(0,t)=u_x(1,t)=  0  & \mbox { in $T$,}
\end{array}
\right.
\end{equation}
we consider the problem
\be\label{liupbneum}
\left\{
\begin{array}{ll}
k_{xx}(x, y, t)-k_{yy}(x, y, t)=(a(y,t)+\l)k(x, y, t),
&\mbox{in $\Omega T$}\\
k_y(x, 0)=0,&\mbox{in $X$}\\
k_x(x, x, t)+k_y(x, x, t)+\frac{d}{d x}(k(x, x, t))=
a(x,t)+\l,&\mbox{in $XT$} \end{array}
\right. \end{equation}
where $\l$ is any constant.  Using the solution $k$,
we then obtain a Dirichlet boundary feedback law
\be\label{dirfdneum}
u(1,t)= - \int_0^1 k(1,y,t)u(y, t)dy \qquad \mbox { in $T$}
\end{equation}
and Neumann boundary feedback law
\be\label{neumfdneum}
u_x(1,t)= -k(1,1)u(1,t)- \int_0^1 k_x(1,y,t)u(y, t)dy
\qquad \mbox { in $T$.}
\end{equation}
With one of the boundary feedback laws, the  system
\begin{equation}
\left\{\begin{array} {ll}
 u_t(x,t)=  u_{xx}(x,t)+a(x,t) u(x,t)
  & \mbox { in $XT,$}\\
u_x(0,t)=0& \mbox { in $T,$}\\
  u(x, 0)=u^0(x)     & \mbox {  in $X$ }
\end{array}
\right.
\label{dcheneum}
\end{equation}
is exponentially stable.  To state this result, we introduce the
compatible conditions for the initial data
\bea\label{compatneu1}
u^0_x(0)&=&0,\quad u^0(1,t)=- \int_0^1 k(1, y,t)u^0(y)dy\\
u^0_x(0)&=&0,\quad u_x^0(1,t)= -k(1,1)u^0(1,t)- \int_0^1 k_x(1,y,t)u^0(y)dy
\label{compatneu2}
\end{eqnarray}

\begin{theorem} \label{cstabthmneum} Assume that  $\l >0$ is any positive
constant and $a(x,t)$ is continuous in
$x$ and analytic in $t$ uniformly in $x$.
For arbitrary  initial  data $u^0(x) \in H^1(0,1) $ with the compatible condition \eqref{compatneu1} or \eqref{compatneu2}, equation \eqref{dcheneum} with either \eqref{dirfdneum} or
 \eqref{neumfdneum} has a unique  solution that satisfies
 \be\nonumber
\|u(t)\|_{H^1} \leq M \|u^0\|_{H^1}e^{- \l  t },\quad t\in (0,t_0)(t_0<1)
\end{equation}
where $M$ is a positive constant independent of $u^0$.
\end{theorem}

\begin{proof}
The proof is the same as that of Theorem~\ref{cstabthm2}.
Only thing we need to do is to show that
problem \eqref{liupbneum} has a unique solution. This is
also considered by Colton \cite{col} and his theorem is quoted as Lemma~\ref{lem3} below.
\end{proof}

\begin{lemma} \label{lem3}  {\bf (Colton)} Suppose that $a(x,t)$ is continuous in
$x$ and analytic in $t$ uniformly in $x$. Then
problem \eqref{liupbneum} has a unique solution which is twice continuously
differentiable in $x,~y$.
\end{lemma}

Similar to Lemma~\ref{lem2}, we have
\begin{lemma} \label{lem4} Let $k(x, y, t)$ be the solution of
problem \eqref{liupbneum} and define the linear bounded operator $K: H^i\! :\!XT_0\to
H^i\! :\!XT_0$ ($i=0,1,2$) by
\be \nonumber
w(x, t)=(Ku)(x, t)=u(x, t)+\int_0^x k(x, y, t)u(y, t)dy.
\end{equation}
Then,
\begin{enumerate}
\item $K$ has a linear bounded inverse $K^{-1}: H^i\! :\!XT_0\to H^i\! :\!XT_0$ ($i=0,1,2$).
\item $K$ converts the system \eqref{dirfdneum} and \eqref{dcheneum} and system \eqref{neumfdneum} and \eqref{dcheneum} into
\begin{equation} \nonumber
\left\{\begin{array} {ll}
 w_t= w_{xx}-\l w
  & \mbox { in $XT,$}\\
  w_x(0,t) = w(1,t)=0& \mbox { in $T,$}\\
  w(x, 0)=w^0(x) & \mbox {  in $X,$ }
\end{array}
\right.
\end{equation}
or
\begin{equation}\nonumber
\left\{\begin{array} {ll}
 w_t= w_{xx}-\l w
  & \mbox { in $XT,$}\\
  w_x(0,t) =  w_x(1,t)=0& \mbox { in $T,$}\\
  w(x, 0)=w^0(x) & \mbox {  in $ X,$ }
\end{array}
\right.
\end{equation}
respectively, where $w^0(x)=u^0(x)+\int_0^x k(x, y, t)u^0(y)dy$.
\end{enumerate}
\end{lemma}

\section{Simulations}\label{simu}
\setcounter {equation} {0}
\setcounter {theorem} {0}
\setcounter {lemma} {0}
To see how the boundary control stabilize the heat equation, we studied the
case of
Dirichlet boundary problem (\ref{linear}) in which \be a(x,t)=x(bt+c),\ee where
$b>0$ and $c>0$ are positive constants. Both the controlled system by the
Dirichlet boundary control (\ref{dirfd}) and uncontrolled
systems are studied. In all the simulation, we use the initial data
\be u^0(x)=10 \( \frac{1}{4}-\(x-\frac{1}{2}\)^2 \)\sin(4\pi
x)+5\(\frac{1}{4}-\(x-\frac{1}{2}\)^2\),\ee which is arbitrarily
chosen. If not stated explicitly, the following values are used in the
simulation:
$\lambda=10,~b=200,~c=5$. The $x$ grid is set to be 100 and the time step size
is set to be $10^{-5}$ but when the figures are drawn, we used larger step
sizes in $x$ and $t$ for clarity.

The corresponding model of the simulation can be viewed as the heat diffusion along a one
dimensional rod, during which heat is also generated unevenly and increasingly. The temperature
of the left end of this rod is fixed zero and another is allowed to adjust to
let the temperature of the rod converge to zero.

\begin{figure}
\mbox{$\begin{array}{cc}
\includegraphics[scale=0.3]{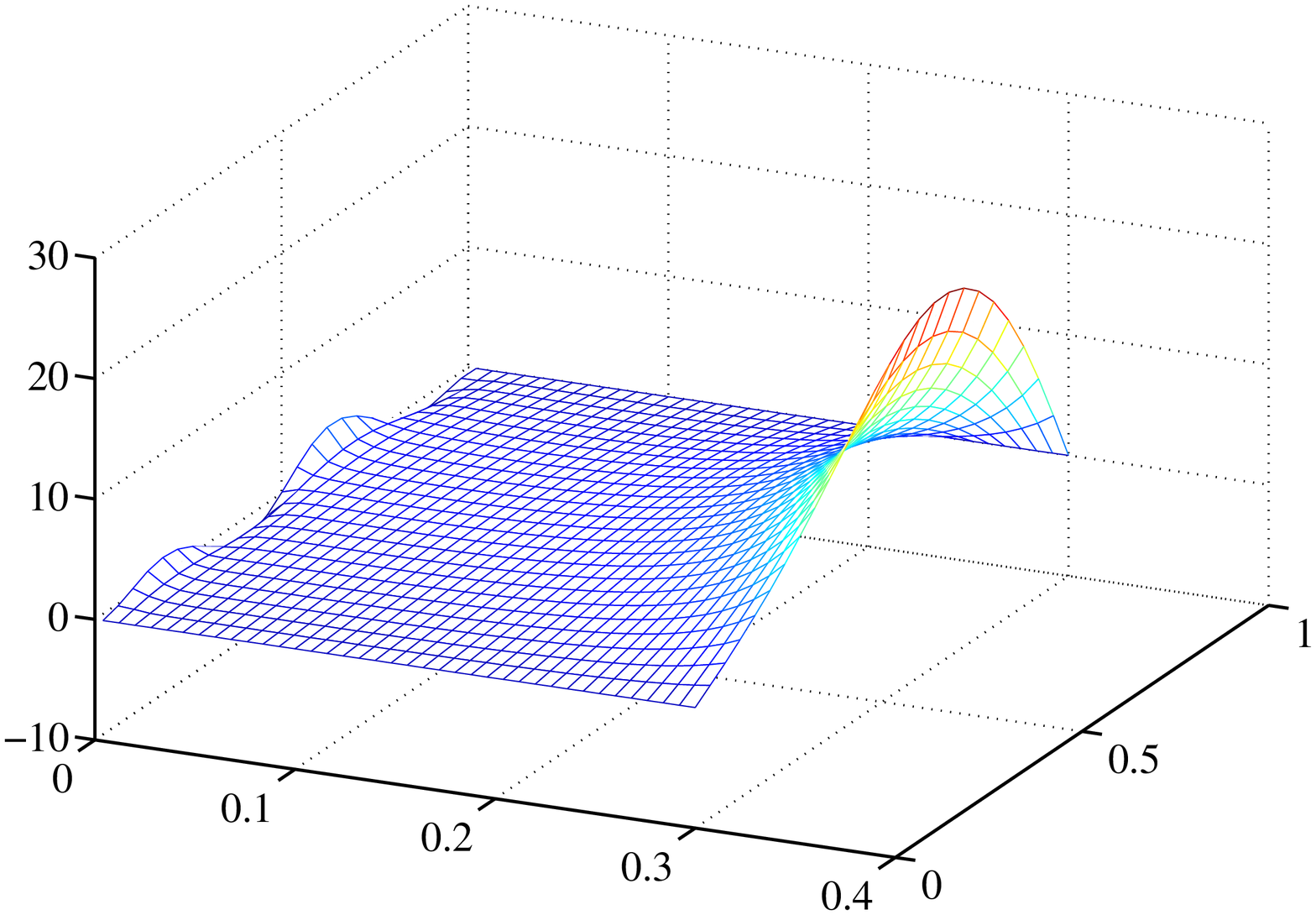} &
\includegraphics[scale=0.3]{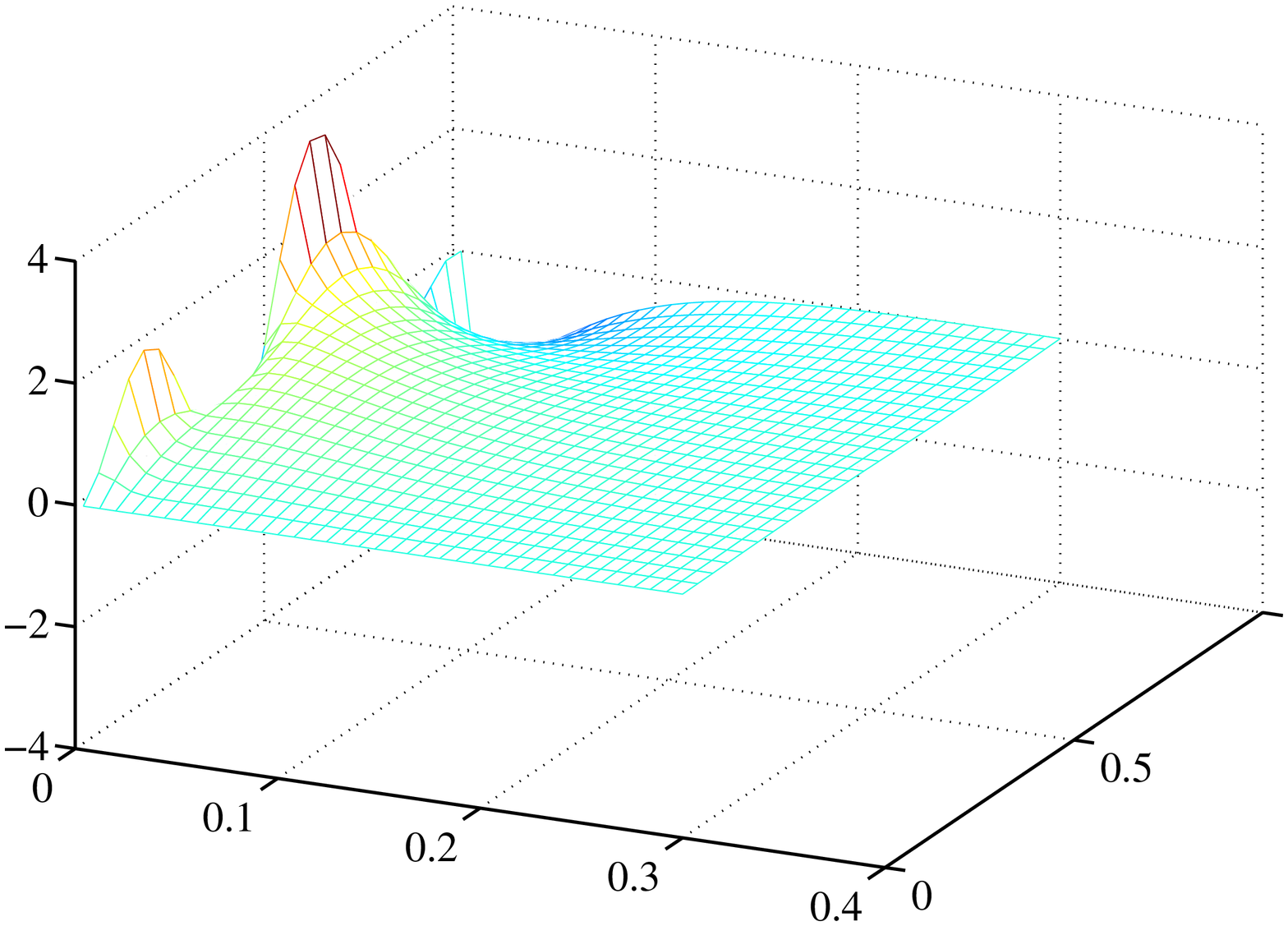}\\
(a)&(b)\\
\includegraphics[scale=0.3]{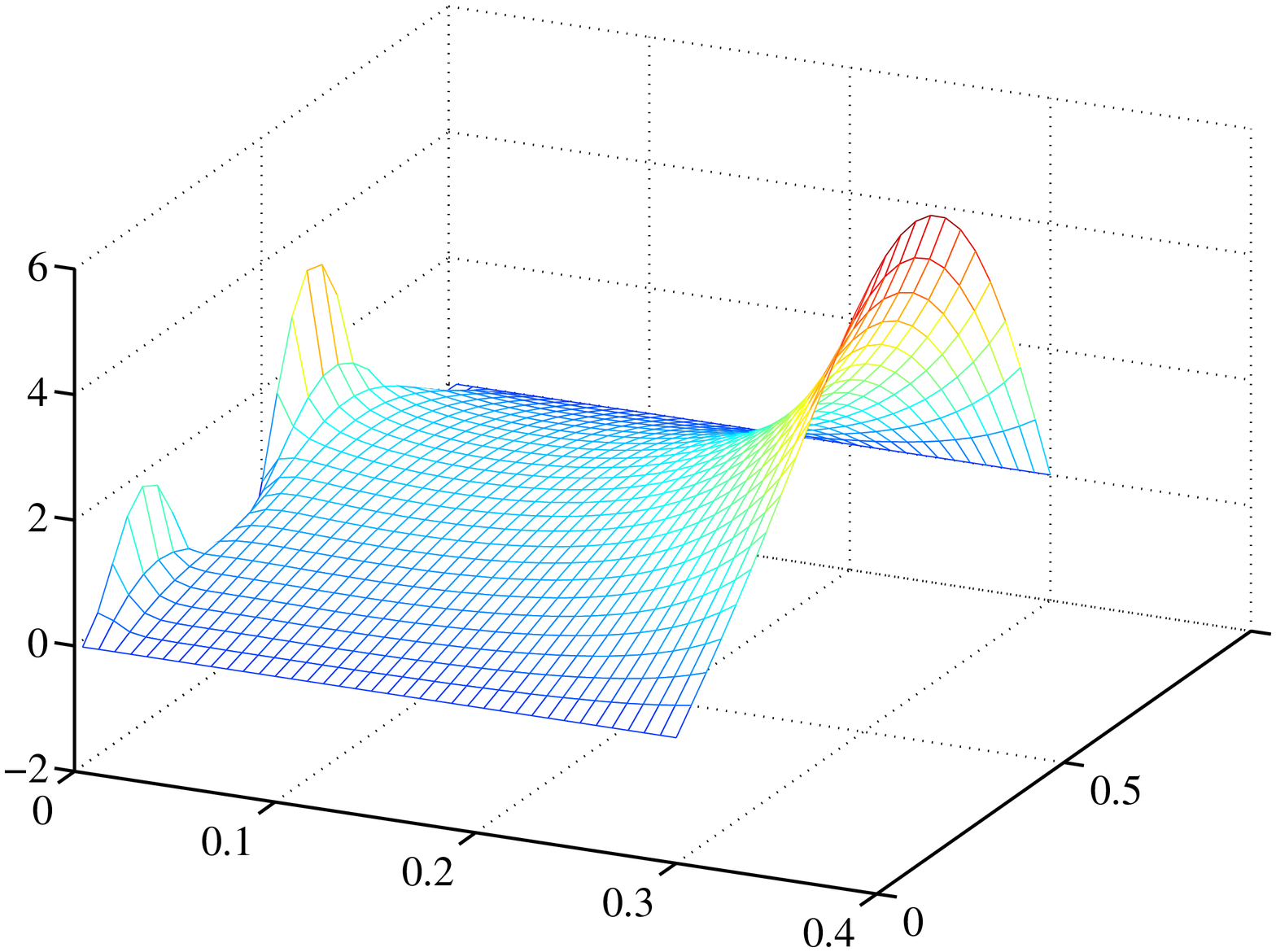} &
\includegraphics[scale=0.3]{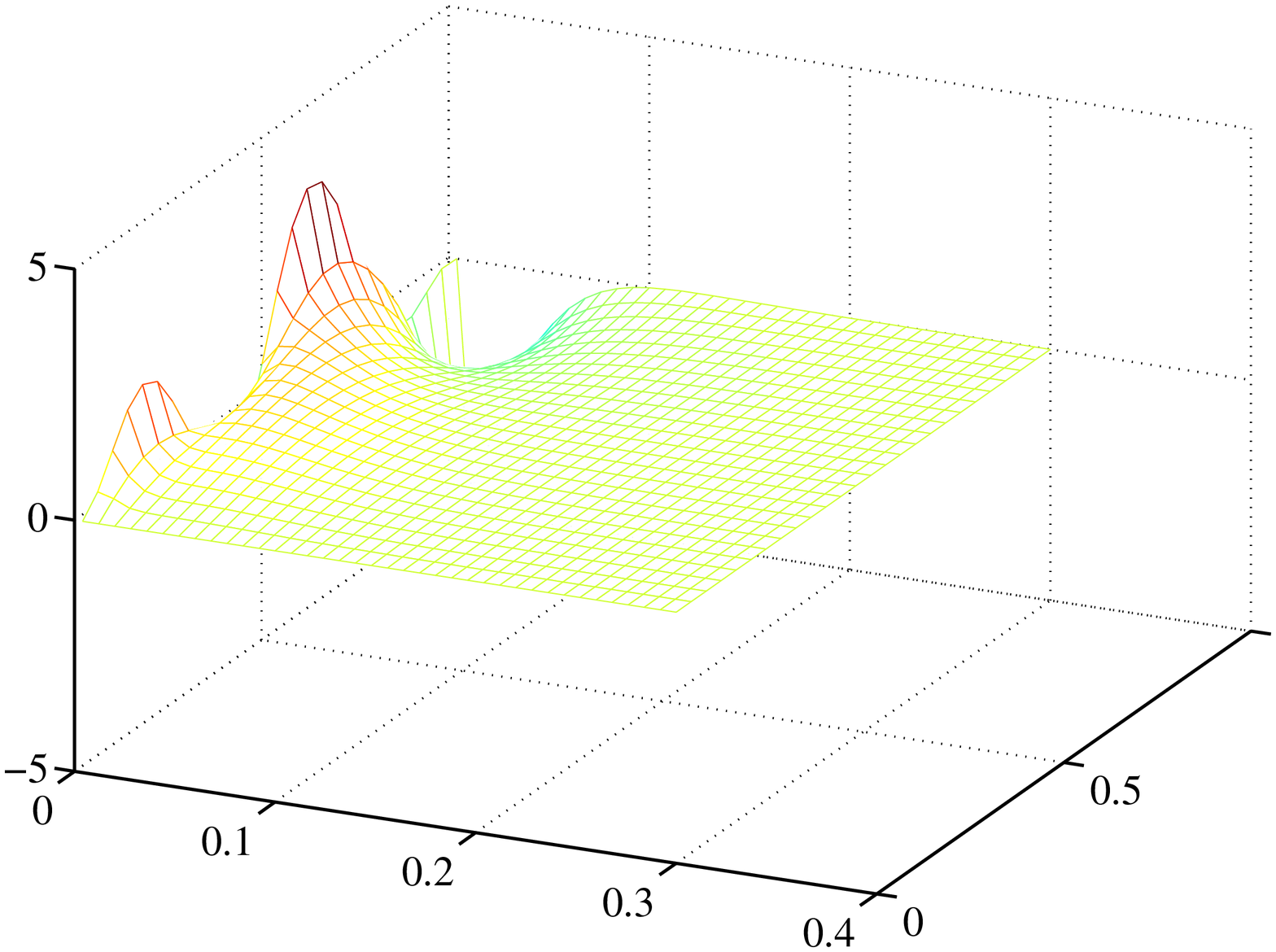}\\
(c)&(d)\\
\includegraphics[scale=0.3]{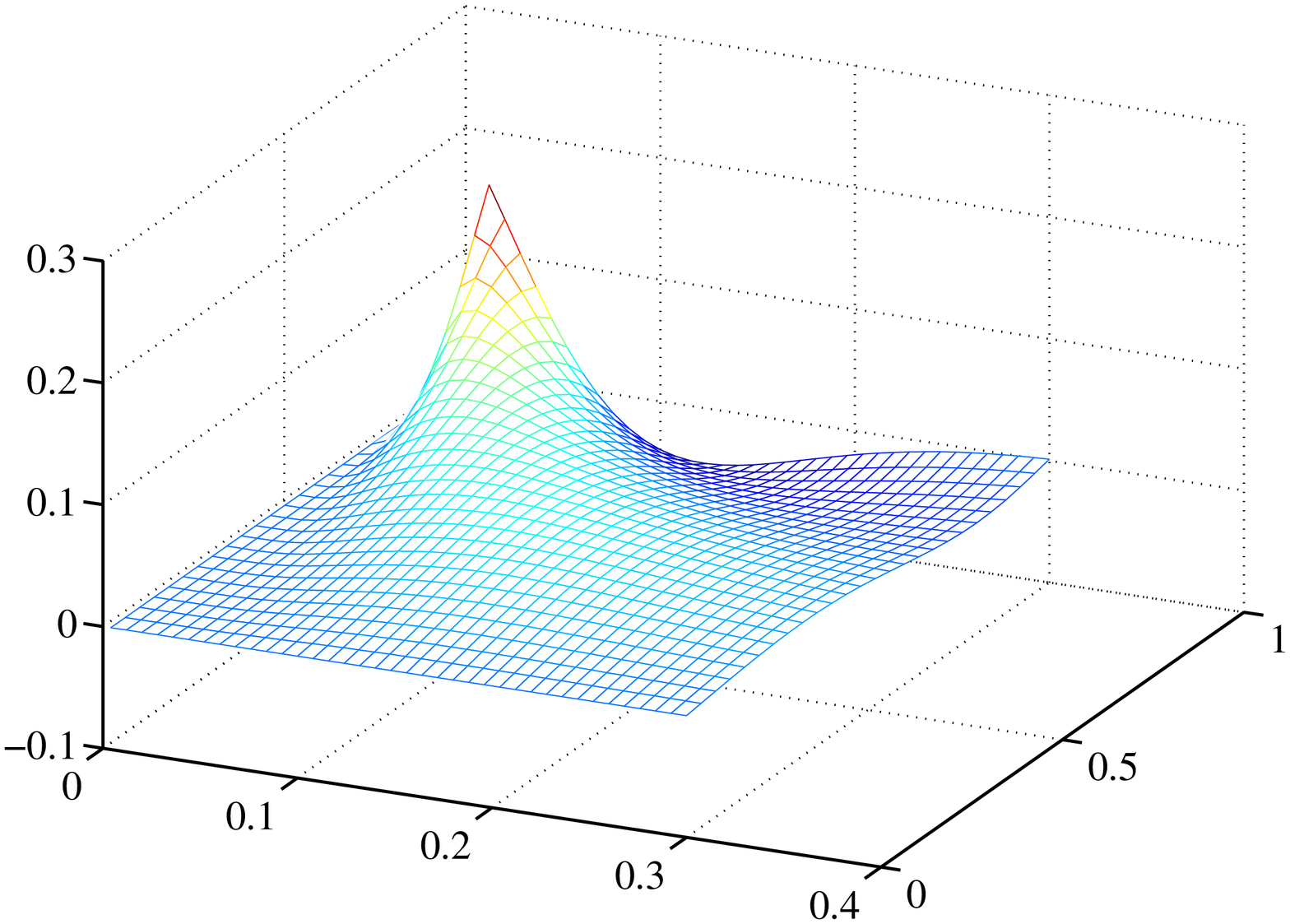} &
\includegraphics[scale=0.3]{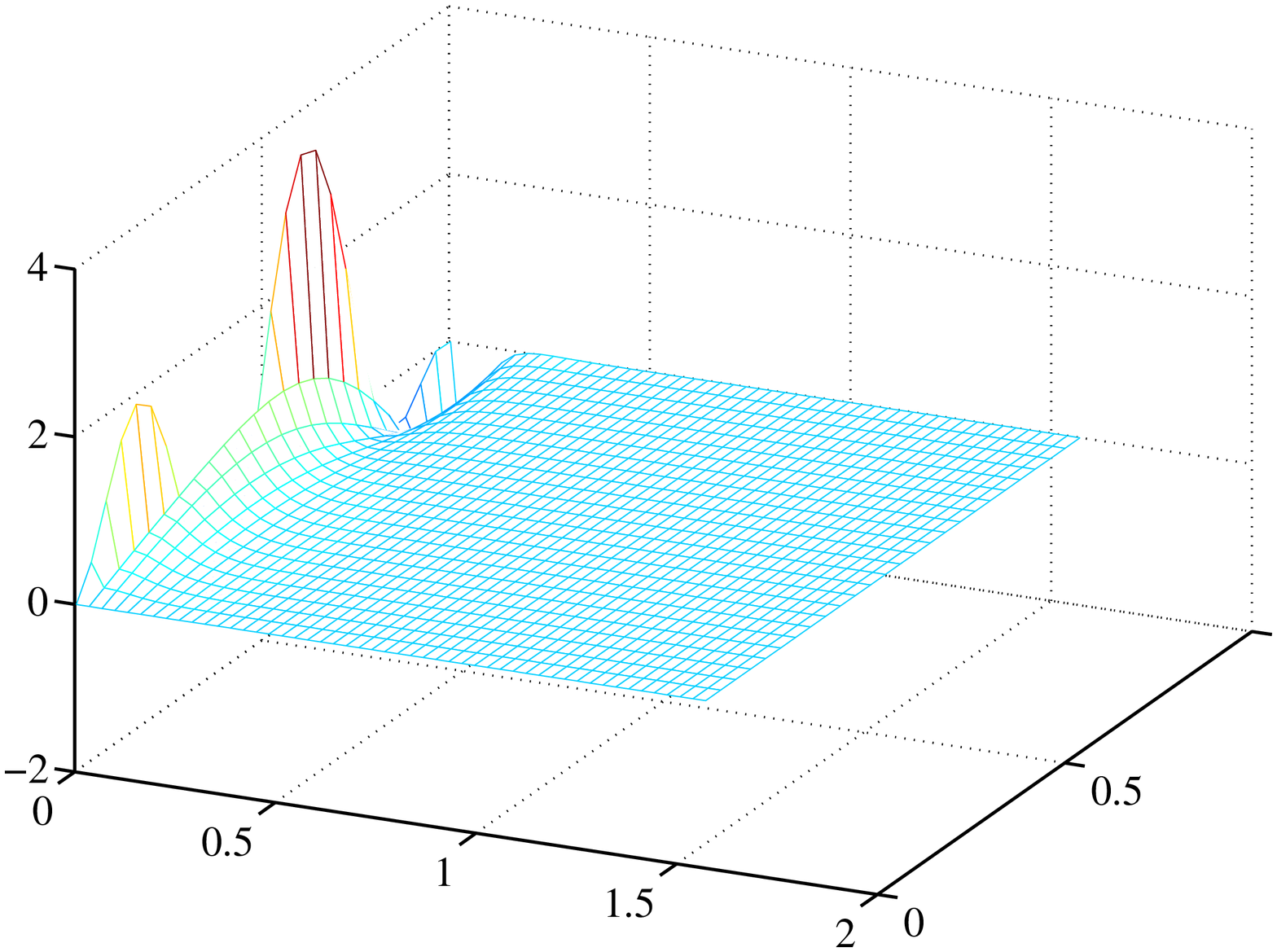}\\
(e)&(f)\\
\end{array}
$} \caption{The temperature of the uncontrolled systems, controlled
systems vs. time; See text for explanation.} \label{figone}
\end{figure}

To reveal the various features and relations of the uncontrolled system,
controlled system and the control inputs, we do the following simulations.
We studied the effects of different $b$'s in $a(x,t)$ with $b=200$ as shown in figure (\ref{figone}.a) and
$b=150$ in figure (\ref{figone}.c) to the uncontrolled system (\ref{linear}). Clearly, (1) the solution to the
uncontrolled system (\ref{linear}) is really divergent with time going and
(2) the larger the coefficient function $a(x,t)$ is (here the bigger the $b$
is), the faster the
solution diverges. The solution to the controlled system with $b=200$ is shown
in Figure (\ref{figone}.b). It is seen that the feedback control can really
stabilize the system. We also studied the effect of different $\lambda$'s in
the kernel $k$, with $\lambda=10$ shown in figure (\ref{figone}.b) and
$\lambda=40$ in figure (\ref{figone}.d). Comparing the $x=1$ edge of these two figures, we can find
that the larger the $\lambda$ is, the faster the solution converges to zero, as
expected from equation (\ref{cstab4}). In all the simulations of
(\ref{figone}.a-d), the kernels are taken to be the first 3 terms in it series
approximation (\ref{series}), that is $G_0+G_1+G_2$. In (\ref{figone}.e), we showed the
difference of solutions to the system controlled by a kernel of first three terms
((\ref{figone}.b)) and the system controlled by a kernel of first four terms.
If denoting these two solutions by $u_3(x,t)$ and $u_4(x,t)$, in figure
(\ref{figone}.e) we showed ($u_3(x,t)-u_4(x,t)$), from which we can see that the high order
approximated kernel can stabilize the system in a faster way. However, this is
not a general conclusion for all $a(x,t)$ since in our case, all terms of
$G_n$ are positive and thus the more terms you use, the better approximation
and better control it gives.

\section{Remarks}
\setcounter {equation} {0}
\setcounter {theorem} {0}
\setcounter {lemma} {0}
As a future work, one can try to stabilize the nonlinear problem
\be u_t(x,t)=u_{xx}(x,t)+a(x,t,u(x,t))u(x,t) \nonumber \ee
where $a$ depends not only on $x,~t$ but also on $u$ itself. However, we
don't know how to prove the
unique existence of the corresponding kernel. Once the unique existence
of the corresponding kernel was shown, all the results in section
(\ref{dir}) and (\ref{neum}) may be generalize to this case.

Another interesting result arises from the simulation. We tried to calculate the
solution to the controlled system (by 3 terms kernel) after $t=1$, which is out
of the time domain in our theorems. From figure (\ref{figone}.f), one can see the system is still stable.
This implies that the condition that
$a(x,t)$ should be analytic in $t$ in theorems (\ref{cstabthm2}) and
(\ref{cstabthmneum}) is only sufficient but not necessary.

\section*{Acknowledgment}
The author thanks Dr. Liu Weijiu for his suggestions to consider
this problem and his valuable discussions. This work was done during
January 2005 and April 2005. The author is grateful to Dr. Vladimir
A. Miranskyy for his support during this work.

\newpage

\end{document}